\documentclass[12pt]{article}
\usepackage{amsmath,amssymb,theorem}
\usepackage{graphicx, enumerate}
\usepackage{chngcntr}
\usepackage{epstopdf}
\usepackage{setspace}

\counterwithin{figure}{section}
\numberwithin{figure}{section}
\counterwithin{table}{section}
\usepackage{hyperref}
\hypersetup{
  colorlinks,
  citecolor=black,
  filecolor=black,
  linkcolor=black,
  urlcolor=black
}

\newtheorem{thm}{Theorem}[section]
\newtheorem{conj}[thm]{Conjecture}

\newtheorem{lem}[thm]{Lemma}

\theorembodyfont{\rmfamily}

\def\pf{\bigskip\noindent {\bf Proof.}~~}

\def\proofsquare{\bigskip\hfill $\blacksquare$}

\newcounter{SectionStep}[section]
\newcommand{\step}[1]{(\refstepcounter{SectionStep}\arabic{SectionStep}) \emph{#1}}
\newcommand{\stepref}[1]{(\ref{#1})}

\newcounter{SideStep}[SectionStep]
\newcommand{\sstep}[1]{(\arabic{SectionStep}.\refstepcounter{SideStep}\arabic{SideStep}) \emph{#1}}
\newcommand{\sstepref}[1]{(\arabic{SectionStep}.\ref{#1})}

\begin{document}
\title{The extremal function for $K_9^=$ minors}
\author{Martin Rolek\thanks{E-mail address: msrolek@wm.edu.} \\
Department of Mathematics\\
College of William \& Mary\\
Williamsburg, VA 23185
}
\date{March 26, 2018}

\maketitle

\begin{abstract}
We prove the extremal function for $K_9^=$ minors, where $K_9^=$ denotes the complete graph $K_9$ with two edges removed.
In particular, we show that any graph with $n$ vertices and at least $6n - 20$ edges either contains a $K_9^=$ minor or is isomorphic to a graph obtained from disjoint copies of $K_8$ and $K_{2, 2, 2, 2, 2}$ by identifying cliques of size 5.
We utilize computer assistance to prove one of our lemmas.
\end{abstract}


\section{Introduction}\label{sec:Intro}

All graphs considered are simple and finite.
We use $V(G)$, $|G|$, $E(G)$, $e(G)$, $\delta(G)$, $\Delta(G)$, and $\overline{G}$ to denote the vertex set, number of vertices, edge set, number of edges, minimum degree, maximum degree, and complement of a graph $G$, respectively.
Given $S, T \subseteq V(G)$, we denote by $e(S, T)$ the number of edges of $G$ with one end in $S$ and one end in $T$.
We denote by $G[S]$ the subgraph of $G$ induced by $S$, and by $G - S$ the subgraph $G[V(G) \setminus S]$ of $G$.
If $S = \{x\}$, we simply write $G - x$ in the latter case.
For $uv \in E(\overline{G})$, we denote by $G + uv$ the graph obtained from $G$ by adding an edge joining $u$ and $v$.
The union (resp. intersection) of two graphs $G$ and $H$, denoted $G \cup H$ (resp. $G \cap H$), is the graph with vertex set $V(G) \cup V(H)$ (resp. $V(G) \cap V(H)$) and edge set $E(G) \cup E(H)$ (resp. $E(G) \cap E(H)$).
The join of two graphs $G$ and $H$, denoted $G \vee H$, is the graph with vertex set $V(G) \cup V(H)$ and edge set $E(G) \cup E(H) \cup \{uv: u \in V(G), v \in V(H)\}$.
If $G$ contains $H$ as a minor, we denote this by $G \ge H$.
$G/xy$ denotes the graph obtained from $G$ by contracting the edge $xy$.
A vertex $x$ is complete to a set $S$ if $x$ is adjacent to every vertex of $S$.
We use $d_G(x)$ to denote the degree of a vertex $x$ in the graph $G$.
Given a subgraph $H$ of $G$, we define $N(H)$ to be the set of vertices of $V(G) \setminus V(H)$ with a neighbor in $V(H)$.

Given a graph property $\mathcal{P}$, the extremal function for $\mathcal{P}$ determines the maximum number of edges a graph on $n$ vertices may have while not satisfying property $\mathcal{P}$.
Extremal graph theory began when Tur\'an~\cite{Turan1941} determined the extremal function for $K_t$ subgraphs.
He further characterized all such graphs attaining this maximum number of edges, the well-known Tur\'an graphs.
Dirac~\cite{Dirac1964a} was the first to consider the natural extension to $K_t$ minors.
When considering $K_t$ minors instead of $K_t$ subgraphs, the problem is much more difficult to solve, and the extremal function is known only for $t \le 9$.
Dirac~\cite{Dirac1964a} showed for $t \le 5$, and Mader~\cite{Mader1968b} for $t \in \{6, 7\}$, that any graph on $n \ge t$ vertices with at least $(t - 2)n - \binom{t - 1}{2} + 1$ edges has a $K_t$ minor.
The case $t = 6$ was also independently shown by Gy\"ori~\cite{Gyori1982}.

For $t \ge 8$, there exist families of graphs with $(t - 2)n - \binom{t - 1}{2} + 1$ edges, but which do not contain a $K_t$ minor.
To describe these families, we define an $(H_1, H_2, k)$-cockade recursively as follows.
Any graph isomorphic to either $H_1$ or $H_2$ is an $(H_1, H_2, k)$-cockade.
Now given two $(H_1, H_2, k)$-cockades $G_1$ and $G_2$, we let $G$ be the graph obtained from $G_1$ and $G_2$ by identifying a clique of size $k$ in $G_1$ with a clique of size $k$ in $G_2$.
Then $G$ is an $(H_1, H_2, k)$-cockade, and every $(H_1, H_2, k)$-cockade can be constructed in this fashion.
If $H_1$ is isomorphic to $H_2$, we simply write $(H_1, k)$-cockade.
J\o rgensen~\cite{Jorgensen1994} showed that any graph on $n \ge 8$ vertices with at least $6n - 20$ edges either has a $K_8$ minor or is a $(K_{2, 2, 2, 2, 2}, 5)$-cockade, and Song and Thomas~\cite{Song2006} showed any graph on $n \ge 9$ vertices with at least $7n - 27$ edges either has a $K_9$ minor or is a $(K_{1, 2, 2, 2, 2, 2}, 6)$-cockade, or is isomorphic to $K_{2, 2, 2, 3, 3}$, settling the cases $t = 8$ and $t = 9$, respectively.
The extremal function for $K_t$ minors remains open for $t \ge 10$.
Note that in a certain sense there is only one minimal counterexample for the case $t = 8$, namely $K_{2, 2, 2, 2, 2}$, and two for the case $t = 9$, namely $K_{1, 2, 2, 2, 2, 2}$ and $K_{2, 2, 2, 3, 3}$.
As pointed out by Song~\cite{Song2004}, there are at least eight minimal counterexamples for the case $t = 10$, and Thomas and Zhu (see~\cite{Thomas2018+}) have conjectured that there are no further minimal counterexamples. 

As a simplification, the extremal function for $K_t^-$ minors has been investigated, where $K_t^-$ is the complete graph on $t$ vertices with one edge removed.
For $t \in \{5, 6\}$, Dirac~\cite{Dirac1964a} showed that any graph on $n \ge t$ vertices with at least $\frac{1}{2}((2t - 5)n - (t - 3)(t - 1))$ edges either has a $K_t^-$ minor or is a $(K_{t - 1}, t - 3)$-cockade.
For larger values of $t$, there is more than one minimal counterexample.
Jakobsen~\cite{Jakobsen1972,Jakobsen1983} showed any graph on $n \ge 7$ vertices with at least $\frac{1}{2}(9n - 24)$ edges either has a $K_7^-$ minor or is a $(K_6, K_{2, 2, 2, 2}, 4)$-cockade, and Song~\cite{Song2005} showed any graph on $n \ge 8$ vertices with at least $\frac{1}{2}(11n - 35)$ edges either has a $K_8^-$ minor or is a $(K_7, K_{1, 2, 2, 2, 2}, 5)$-cockade.
The extremal function for $K_t^-$minors remains open for $t \ge 9$.

In this paper, we will consider $K_t^=$ minors, where $K_t^=$ denotes the complete graph on $t$ vertices with two edges removed.
Note that there are two nonisomorphic graphs $K_t^=$, depending on whether the removed edges share a common end or not.
Let $\mathcal{K}_t^=$ denote the family consisting of the two nonisomorphic graphs $K_t^=$.
Throughout this paper, we will use the following conventions.
We say a graph $G$ has no $K_t^=$ minor if $G$ does not contain $K$ as a minor for any $K \in \mathcal{K}_t^=$, and we say that $G$ has a $K_t^=$ minor if $G$ contains $K$ as a minor for some $K \in \mathcal{K}_t^=$.
Dirac~\cite{Dirac1964a} proved the following for $t \in \{5, 6\}$, and Jakobsen~\cite{Jakobsen1971,Jakobsen1972} proved the cases $t \in \{7, 8\}$.

\begin{thm}\label{thm:K8=Extremal}
(Dirac~\cite{Dirac1964a}, Jakobsen~\cite{Jakobsen1971,Jakobsen1972})
For $t \in \{5, \dots, 8\}$, if $G$ is a graph with $|G| \ge t - 1$ and at least $(t - 3)n - \frac{1}{2}(t - 1)(t - 4)$ edges, then either $G \ge K_t^=$ or $G$ is a $(K_{t - 1}, t - 4)$-cockade.
\end{thm}

Our main result is to extend Theorem~\ref{thm:K8=Extremal} to the case $t = 9$ as follows.

\begin{thm}\label{thm:K9=Extremal}
If $G$ is a graph with $|G| \ge 8$ and at least $6|G| - 20$ edges, then either $G \ge K_9^=$ or $G$ is a $(K_8, K_{2, 2, 2, 2, 2}, 5)$-cockade.
\end{thm}

Note that for this case there are now two minimal counterexamples to consider.
In Section~\ref{sec:Prelim} we prove several results necessary for the proof of Theorem~\ref{thm:K9=Extremal}, which we present in Section~\ref{sec:K9=Extremal}.

Our primary motivation for studying the extremal functions for $K_t$, $K_t^-$, and $K_t^=$ minors is their integral use in proving results related to Hadwiger's conjecture~\cite{Hadwiger1943}, which claims that every graph with no $K_t$ minor is $(t - 1)$-colorable.
Hadwiger's conjecture is easily true for $t \le 3$.
The case $t = 4$ was shown by both Hadwiger~\cite{Hadwiger1943} and Dirac~\cite{Dirac1952}, and a short proof was given much later by Woodall~\cite{Woodall1992}.
For $t = 5$, Wagner~\cite{Wagner1937} showed that Hadwiger's conjecture is equivalent to the Four Color Theorem, and for $t = 6$, Robertson, Seymour, and Thomas~\cite{Robertson1993} showed the same.
The conjecture remains open for $t \ge 7$, although there are some partial results as follows.
Albar and Gon\c calves~\cite{Albar2017} showed for $t \in \{7, 8\}$, and the present author and Song~\cite{Rolek2017} for $t = 9$, that every graph with no $K_t$ minor is $(2t - 6)$-colorable.
An alternate proof for the cases $t \in \{7, 8\}$ is also provided in~\cite{Rolek2017}.
By noticing that known minimal counterexamples to the extremal function for $K_t$ minors for $t \ge 7$ are all $(t - 1)$-colorable, the present author and Song~\cite{Rolek2017} proved the first general result on coloring graphs with no $K_t$ minor for all $t \ge 6$, provided a suitable conjecture holds as follows.

\begin{conj}\label{conj:ExtremalColor}
(Rolek and Song~\cite{Rolek2017})
For every $t \ge 1$, every graph on $n$ vertices with at least $(t - 2)n - \binom{t - 1}{2} + 1$ edges either has a $K_t$ minor or is $(t - 1)$-colorable.
\end{conj}

\begin{thm}\label{thm:KtColoring}
(Rolek and Song~\cite{Rolek2017})
For $t \ge 6$, if Conjecture~\ref{conj:ExtremalColor} is true, then every graph with no $K_t$ minor is $(2t - 6)$-colorable.
\end{thm}

The chromatic number of graphs without $K_t^-$ minors and $K_t^=$ minors has also been investigated.
Jakobsen~\cite{Jakobsen1972,Jakobsen1983} showed for $t = 7$, and the present author and Song~\cite{Rolek2017} for $t = 8$, that if $G$ has no $K_t^-$ minor, then $G$ is $(2t - 7)$-colorable, and if $G$ has no $K_t^=$ minor, then $G$ is $(2t - 8)$-colorable.
Most recently, the present author~\cite{Rolek2018b} has used Theorem~\ref{thm:K9=Extremal} to show the following.

\begin{thm}\label{thm:K9=Coloring}
(Rolek~\cite{Rolek2018b})
If $G$ has no $K_9^=$ minor, then $G$ is 10-colorable.
\end{thm}

\section{Preliminaries}\label{sec:Prelim}

We begin this section with four results on $(K_8, K_{2, 2, 2, 2, 2}, 5)$-cockades which will be useful early in Section~\ref{sec:K9=Extremal}.

\begin{lem}\label{lem:K85Cockade+e}
Let $G$ be a $(K_8, K_{2, 2, 2, 2, 2}, 5)$-cockade, and let $x$ and $y$ be nonadjacent vertices in $G$.
Then $G + xy \ge K_9^=$.
\end{lem}

\pf
We proceed by induction on $|G|$.
It is easy to see that the statement holds if $G$ is isomorphic to $K_{2, 2, 2, 2, 2}$.
Hence we may assume that $G$ is obtained from $H_1$ and $H_2$ by identifying a $K_5$, where both $H_1$ and $H_2$ are $(K_8, K_{2, 2, 2, 2, 2}, 5)$-cockades.
If both $x, y \in V(H_i)$ for some $i \in \{1, 2\}$, then $G \ge K_9^=$ by induction.
Thus we may assume that $x \in V(H_1) \setminus V(H_2)$ and $y \in V(H_2) \setminus V(H_1)$.
If there exists $z \in V(H_1 \cap H_2)$ such that $xz \notin E(G)$, then by contracting the component of $H_2 - V(H_1 \cap H_2)$ containing $y$ onto $z$ and deleting all other components of $H_2 - V(H_1 \cap H_2)$, we see that $G \ge H_1 + xz$, and the resulting graph contains a $K_9^=$ minor by induction.
Hence we may assume that $x$ is complete to $V(H_1 \cap H_2)$, and similarly that $y$ is complete to $V(H_1 \cap H_2)$.
Since $G[(V(H_1 \cap H_2)) \cup \{x\}]$ is isomorphic to $K_6$, it follows that there is a $K_8$-subgraph $H'$ of $H_1$ such that $(V(H_1 \cap H_2)) \cup \{x\} \subseteq V(H')$.
Then $G[V(H') \cup \{y\}]$ is a $K_9^=$-subgraph in $G + xy$.
\proofsquare

From Lemma~\ref{lem:K85Cockade+e} we get the following.

\begin{lem}\label{lem:K85Cockade+v}
If $G'$ is a $(K_8, K_{2, 2, 2, 2, 2}, 5)$-cockade, and $G$ is the graph obtained from $G'$ by adding a new vertex joined to at least six vertices of $G'$, then $G \ge K_9^=$.
\end{lem}

\pf
Let $G$ and $G'$ be as in the statement, and say $v \in V(G) \setminus V(G')$.
If there exist $x, y \in N_{G}(v)$ such that $xy \notin E(G)$, then $G \ge G' + xy \ge K_9^=$ by Lemma~\ref{lem:K85Cockade+e}.
Hence $N_{G}(v)$ is complete, and so there exists a $K_8$-subgraph $H$ of $G'$ such that $N_{G}(v) \subseteq V(H)$.
But then $G[V(H) \cup \{v\}]$ is isomorphic to $K_9^=$.
\proofsquare

\begin{lem}\label{lem:K85CockadeContract}
Let $G$ be a graph with $\delta(G) \ge 7$.
Let $xy \in E(G)$ such that $x$ and $y$ have at least 5 common neighbors.
If $G/xy$ is a $(K_8, K_{2, 2, 2, 2, 2}, 5)$-cockade, then $G \ge K_9^=$.
\end{lem}

\pf
We proceed by induction on $|G|$.
The statement is easy to verify if $G/xy$ is isomorphic to $K_8$.
Assume $G/xy$ is isomorphic to $K_{2, 2, 2, 2, 2}$.
Say $V(G/xy) = \{v_1, w_1, \dots, v_5, w_5\}$, where $v_i w_i \notin E(G/xy)$ for $i \in \{1, \dots, 5\}$, and $w_5$ is the vertex obtained by contracting the edge $xy$ of $G$.
We may assume by Pigeonhole that $v_1, w_1, v_2, v_3$ are common neighbors of $x$ and $y$ in $G$.
If $w_2$ is also a common neighbor of $x$ and $y$, then we may assume by symmetry that $xw_3, yw_4 \in E(G)$ since $\delta(G) \ge 7$ and each of $w_3, w_4$ must be adjacent to at least one of $x, y$.
Now by contracting the edges $v_1 v_5$ and $w_3 w_4$, and noting that $v_4$ also must be adjacent to at least one of $x, y$, we see that $G \ge K_9^=$.
Hence we may assume $w_2$ is not a common neighbor of $x$ and $y$, and by symmetry neither is $w_3$.
Then $v_4$, say, is a common neighbor of $x$ and $y$.
Now since $\delta(G) \ge 7$ and each of $w_2, w_3, w_4$ must be adjacent to at least one of $x, y$, we may assume that $xw_2, yw_3 \in E(G)$.
By contracting the edges $v_1 v_5$ and $w_2 w_3$, we again see $G \ge K_9^=$.

Therefore we may assume that $G/xy$ is a $(K_8, K_{2, 2, 2, 2, 2}, 5)$-cockade obtained from $H_1$ and $H_2$ by identifying a $K_5$, where both $H_1$ and $H_2$ are $(K_8, K_{2, 2, 2, 2, 2}, 5)$-cockades.
Say $v \in V(G/xy)$ is the vertex resulting from contracting the edge $xy$ of $G$.
Let $S = V(H_1 \cap H_2)$.
Let $H_i'$ be the subgraph of $G$ induced by $(V(H_i) \setminus \{v\}) \cup \{x, y\}$ for $i \in \{1, 2\}$.
If $v \in V(H_1) \setminus V(H_2)$, then every common neighbor of $x$ and $y$ in $G$ belongs to $V(H_1')$.
Further, every vertex of $S$ has at least 7 neighbors in $V(H_1')$.
Hence $G \ge K_9^=$ by induction applied to $H_1'$.
Thus $v \notin V(H_1) \setminus V(H_2)$ and similarly, $v \notin V(H_2) \setminus V(H_1)$.
Therefore $v \in S$.
Let $Z$ denote the set of common neighbors of $x$ and $y$ in $G$.
We may assume that $|Z \cap V(H_1)| \ge |Z \cap V(H_2)|$, and in particular that $Z \cap (V(H_1) \setminus V(H_2)) \ne \emptyset$ since $|Z| \ge 5$ and $|S \setminus \{v\}| = 4$.
We may also assume that $d_{H_1'}(x) \ge d_{H_1'}(y)$.
Let $Z'$ be a set of neighbors of $y$ in $H_2' - x$ chosen maximal subject to $|Z'| \le 4$ and then further such that $|Z \cap Z'|$ is maximum.
Say $|Z'| = r$.
Since $G/xy$ is a $(K_8, K_{2, 2, 2, 2, 2}, 5)$-cockade, it is 5-connected.
Thus there exist $r$ disjoint paths $P_1, \dots, P_r$ in $G/xy - v$ with one end in $Z'$ and the other end in $S \setminus \{v\}$.
Then each path $P_i$ is also a path in $G$.
Let $H^*$ be the graph obtained from $G$ by contracting each path $P_i$ onto its end in $S$ and deleting all vertices of $V(H_2) \setminus V(H_1 \cup P_1 \cup \dots \cup P_r)$.
Then by the choice of $Z'$, $x$ and $y$ have at least 5 common neighbors in $H^*$.
Also, $H^* \ge H_1$.
If $d_{H^*}(y) \ge 7$, then $\delta(H^*) \ge 7$, so $G \ge H^* \ge K_9^=$ by induction applied to $H^*$.
Otherwise, $d_{H^*}(y) = 6$, $y$ has only one neighbor in $V(H^*) \setminus S$, say $w \in Z$, and every vertex of $S \setminus \{v\}$ is a common neighbor of $x$ and $y$.
Since $H_1$ is a cockade and $H_1[S]$ is isomorphic to $K_5$, some subgraph $K$ of $H_1$ contains $S$, where $K$ is either isomorphic to $K_8$ or $K_{2, 2, 2, 2, 2}$.
Then there exists a path $Q$ in $H_1$ with one end $w$ and the other end in $V(K)$ such that $Q$ avoids $S$, and $Q$ has no internal vertices in $V(K)$ (possibly $Q$ consists only of the vertex $w$).
By contracting $Q$ onto its end in $V(K)$, we may assume that $y$ has one neighbor in $V(K) \setminus S$.
Let $K^*$ be the subgraph of $H^*$ induced by $V(K) \cup \{x, y\} \setminus \{v\}$.
Note that every vertex of $V(K)$ is adjacent to $x$, except the nonneighbor of $v$ in $K$ if $K$ is isomorphic to $K_{2, 2, 2, 2, 2}$, because each such vertex is adjacent to at least one of $x$ or $y$, and every neighbor of $y$ in $V(K)$ is a common neighbor of $x$ and $y$.
Hence, $K^* - y$ is isomorphic to $K$, that is $K^* - y$ is isomorphic to $K_8$ or $K_{2, 2, 2, 2, 2}$ and in particular is a $(K_8, K_{2, 2, 2, 2, 2}, 5)$-cockade.
Since $y$ has six neighbors in $V(K^*)$, it follows from Lemma~\ref{lem:K85Cockade+v} that $G \ge K^* \ge K_9^=$.
\proofsquare

It is easy to verify the following, so the details are omitted.

\begin{lem}\label{lem:K85CockadeEdges}
If $G$ is a $(K_8, K_{2, 2, 2, 2, 2}, 5)$-cockade, then $e(G) = 6|G| - 20$.
\end{lem}

Given a graph $G$ and a set $T = \{v_1, \dots, v_4\} \subseteq{V(G)}$, we say that $G$ has a $K_4$ minor rooted at $T$ if, for $i \in \{1, \dots, 4\}$, $G$ has disjoint subsets $T_i \subseteq V(G)$ such that $G[T_i]$ is connected, $T_i \cap T = \{v_i\}$, and there exist some vertices $u_i \in T_i$ and $u_j \in T_j$ such that $u_i  u_j \in E(G)$ for all $j \in \{1, \dots, 4\}$ with $j \ne i$.
The concept of rooted $K_4$ minors is an extension of earlier work by Seymour~\cite{Seymour1980} and Thomassen~\cite{Thomassen1980} on 2-linked graphs.
Robertson, Seymour, and Thomas utilized rooted $K_4$ minors in their proof of Hadwiger's conjecture for graphs with no $K_6$ minor~\cite{Robertson1993}, and the next result is a simplified restatement of (2.6) from that same paper.
For a complete characterization of graphs with rooted $K_4$ minors, see Fabila-Monroy and Wood~\cite{Fabila-Monroy2013}.

\begin{thm}\label{thm:RootedK4} (Robertson, Seymour, and Thomas~\cite{Robertson1993})
If $G$ is 4-connected and $T \subseteq V(G)$ with $|T| = 4$, then either $G$ has a $K_4$-minor rooted at $T$, or $G$ is planar with $e(G) \le 3|G| - 7$.
\end{thm}

\begin{lem}\label{lem:11VK5Subgraph}
Let $G$ be a graph with $|G| \in \{7, \dots, 11\}$ and $\delta(G) \ge 6$.
If $G$ is not 4-connected, then $G$ contains $K_5$ as a subgraph.
If $G$ is not 5-connected, then $G$ contains $K_5^-$ as a subgraph.
\end{lem}

\pf
The statement is clearly true if $G$ is a complete graph, so we may assume that $|G| \in \{8, \dots, 11\}$.
Let $S$ be a minimum separating set in $G$, and let $G_1, G_2$ be proper subgraphs of $G$ such that $G_1 \cup G_2 = G$ and $G_1 \cap G_2 = G[S]$.
Since $\delta(G) \ge 6$, it follows that $|G_i| \ge 7$ for $i \in \{1, 2\}$.
Since $|G| = |G_1| + |G_2| - |S| \ge 14 - |S|$, we have $|S| \ge 5$ for $|G| \in \{8, 9\}$, and $|S| \ge 4$ for $|G| = 10$.
Assume first that $G$ is not 4-connected.
Then $|G| = 11$ and $|S| = 3$, so $|G_1| = 7$.
Any vertex $u \in V(G_1)\setminus S$ is complete to $V(G_1) \setminus \{u\}$.
Since $|G_1 - S| = 4$, the subgraph of $G$ induced by $(V(G_1) \setminus S) \cup \{v\}$ for any $v \in S$ is isomorphic to $K_5$.
Now assume that $G$ is 4-connected, but not 5-connected.
Then $|G| \in \{10, 11\}$.
At least one of $G_1$, $G_2$, say $G_1$, must satisfy $|G_1| = 7$.
Any vertex $u \in V(G_1)\setminus S$ is complete to $V(G_1) \setminus \{u\}$.
Since $|G_1 - S| \ge 3$, the subgraph of $G$ induced by $(V(G_1) \setminus S) \cup \{v_1, v_2\}$ for any $v_1, v_2 \in S$ contains $K_5^-$ as a subgraph.
\proofsquare

\begin{lem}\label{lem:11VK5}
Let $G$ be a graph with $|G| \in \{7, \dots, 11\}$ and $\delta(G) \ge 6$.
Then either $G$ contains $K_5$ as a subgraph, or for any set $T \subseteq V(G)$ with $|T| = 5$ there exists $v \in T$ such that $G - v$ has a $K_4$ minor rooted at $T \setminus \{v\}$.
\end{lem}

\pf
By Lemma~\ref{lem:11VK5Subgraph}, if $G$ is not 4-connected, then $G$ contains $K_5$ as a subgraph and we are done.
Thus we may assume $G$ is 4-connected.
Let $T \subseteq V(G)$ such that $|T| = 5$.
Assume first that $G$ is 5-connected, and let $v \in T$ be arbitrary.
Then $G - v$ is a $4$-connected graph on $|G| - 1$ vertices with $\delta(G) \ge 5$.
Note that $e(G - v) \ge \left\lceil\frac{1}{2}(|G| - 1)\delta(G)\right\rceil \ge 3(|G| - 1) - 6$, and it follows from Theorem~\ref{thm:RootedK4} that $G - v$ has a $K_4$ minor rooted at $T \setminus \{v\}$.
Therefore we may assume $G$ is not 5-connected, and so by Lemma~\ref{lem:11VK5Subgraph} $G$ contains $K_5^-$ as a subgraph.
In particular, $G$ contains a $K_4$ subgraph $H$.
Since $G$ is 4-connected, there exist four disjoint paths $P_1, P_2, P_3, P_4$ with one end in $V(H)$, one end in $T$, and no internal vertices in $V(H) \cup T$.
Let $v \in T$ be the unique vertex not met by any $P_i$.
Then contracting each path $P_i$ to a single vertex gives a $K_4$ minor rooted at $T$ in $G - v$.
\proofsquare

\begin{lem}\label{lem:EdgeorK14}
Let $G$ be a graph with $|G| \in \{8, \dots, 11\}$ and $\delta(G) \ge 6$, and let $v_1, \dots, v_6 \in V(G)$ be distinct.
If $v_1 v_2 \notin E(G)$, then there exists some component $C$ of $G - \{v_1, \dots, v_6\}$ such that either $\{v_1, v_2\} \subseteq N_G(C)$ or $\{v_3, \dots, v_6\} \subseteq N_G(C)$.
\end{lem}

\pf
Suppose for a contradiction that no such component exists.
Put $T = \{v_1, \dots, v_6\}$.
Since $\delta(G) \ge 6$ and $|T| = 6$, every vertex of $T$ has at least one neighbor in $V(G) \setminus T$.
Thus it must be the case that $G - T$ is disconnected.
Let $C_1, C_2$ be distinct components of $G - T$.
Then for $i \in \{1, 2\}$, any vertex of $C_i$ is adjacent to at most four vertices of $T$, and it follows that $|V(C_i)| \ge 3$ since $\delta(G) \ge 6$.
But then $|G| \ge |T| + |C_1| + |C_2| \ge 12$, a contradiction.
\proofsquare

The following lemma is proved by computer search.
The details of this search can be found in~\cite{Rolek2018c}.

\begin{lem}\label{lem:K7=UK1Computer}
If $G$ is a graph with $|G| \in \{9, 10, 11\}$ and $\delta(G) \ge 6$, then either $G$ contains $K_7^= \cup K_1$ as a minor, or $G$ is isomorphic to one of the five graphs $\overline{C_5} \vee \overline{C_4}$, $\overline{C_9}$, $K_{3, 3, 3}$, $\overline{C_6} \vee \overline{K_3}$, or $\overline{P}$, where $\overline{P}$ is the complement of the Petersen graph $P$.
Furthermore, the graphs $\overline{C_5} \vee \overline{C_4}$ and $\overline{C_9}$ contain $K_7^-$ as a minor, and the graphs $K_{3, 3, 3}$, $\overline{C_6} \vee \overline{K_3}$, and $\overline{P}$ are all edge maximal with respect to not having a $K_7^= \cup K_1$ minor.
\end{lem}

\begin{figure}[t]\centering
\includegraphics[width=125px]{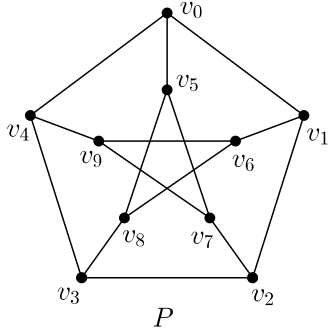}
\caption{\label{fig:Petersen} The Petersen graph $P$.}
\end{figure}

The remaining results in this section will all be used in the proof of Theorem~\ref{thm:K9=Extremal} to handle the counterexamples $K_{3, 3, 3}$, $\overline{C_6} \vee \overline{K_3}$, and $\overline{P}$ of Lemma~\ref{lem:K7=UK1Computer}.

\begin{lem}\label{lem:K333}
Let $G$ be isomorphic to either $K_{3, 3, 3}$ or $\overline{C_6} \vee \overline{K_3}$.
If $e_1, e_2, e_3$ are distinct missing edges of $G$ such that $e_1, e_2$ share a common end, then $G + \{e_1, e_2, e_3\} \ge K_8^=$.
\end{lem}

\pf
Since $e_1, e_2, e_3$ are not independent, there exist $u, v \in V(G)$ such that $uv \in E(G)$ and $u$, $v$ are both incident to two missing edges of $G + \{e_1, e_2, e_3\}$.
By contracting $uv$ we obtain a $K_8^=$ minor of $G + \{e_1, e_2, e_3\}$.
\proofsquare

\begin{lem}\label{lem:Petersen}
Let $e_1, e_2, e_3$ be three distinct missing edges of $\overline{P}$ such that no vertex of $\overline{P}$ is incident to every $e_i$.
If either the $e_i$ do not all belong to the same 5 vertex cycle in $P$ or the $e_i$ induce a 4 vertex path in $P$, then $\overline{P} + \{e_1, e_2, e_3\} \ge K_8^=$.
Furthermore, the graph obtained from $\overline{P}$ by adding any four missing edges contains $K_8^=$ as a minor.
\end{lem}

\pf
Assume first that $e_1, e_2, e_3$ all belong to a 4 vertex path in $P$.
Without loss of generality, $\{e_1, e_2, e_3\} = \{v_0 v_1, v_1 v_2, v_2 v_3\}$, where the vertices of $P$ are labeled as in Figure~\ref{fig:Petersen}.
By contracting the edges $v_0 v_6$ and $v_3 v_7$, we see that $\overline{P} + \{v_0 v_1, v_1 v_2, v_2 v_3\} \ge K_8^=$.
So we may now assume that $e_1, e_2, e_3$ do not all belong to some 5 vertex cycle in $P$.
If the edges $e_i$ are pairwise disjoint, then
without loss of generality, $\{e_1, e_2, e_3\}$ is one of the sets $\{v_0 v_5, v_2 v_3, v_6 v_9\}$, $\{v_0 v_1, v_2 v_3, v_5 v_8\}$, $\{v_0 v_5, v_1 v_2, v_3 v_4\}$, or $\{v_0 v_5, v_1 v_6, v_3 v_4\}$.
In each case, it is straightforward to show that $\overline{P} + \{e_1, e_2, e_3\} \ge K_8^=$.
So we may assume that two $e_i$ are incident to the same vertex, say $\{e_1, e_2\} = \{v_0 v_1, v_0 v_4\}$.
Then $e_3$ is either incident to $v_5$ or not, and by symmetry we may assume $e_3$ is either $v_5 v_7$ or $v_6 v_8$.
Both of these cases are also straightforward to verify.
For the second part of the statement, it is easy to see that given any four missing edges of $\overline{P}$, some three of those edges satisfy the conditions of the first part of the statement.
\proofsquare

\begin{lem}\label{lem:Petersen7}
Let $A_1, A_2 \subseteq V(\overline{P})$ such that $|A_1| \ge |A_2| \ge 7$.
Then for $i \in \{1, 2\}$ there exist vertices $v_i \in A_i$ such that the graph obtained from $\overline{P}$ by adding all missing edges incident to $v_i$ in $\overline{P}[A_i]$ contains $K_8^=$ as a minor.
\end{lem}

\pf
Let the vertices of $P$ be as labeled in Figure~\ref{fig:Petersen}.
Since $|A_1| \ge 7$, there must exist two missing edges of $\overline{P}[A_1]$ which share a common end, say $\{v_0, v_1, v_4\} \subseteq A_1$.
From Lemma~\ref{lem:Petersen}, it follows that if $\overline{P}[A_2]$ contains any missing edge $e$ other than $v_0 v_1$, $v_0 v_4$, $v_0 v_5$, $v_2 v_3$, or $v_6 v_9$, then we are done by adding all missing edges in $\overline{P}[A_1]$ incident to $v_0$ and all missing edges in $\overline{P}[A_2]$ incident to one end of $e$.
Since $|A_2| \ge 7$, such an edge $e$ must exist.
\proofsquare

\section{Proof of Theorem~\ref{thm:K9=Extremal}}\label{sec:K9=Extremal}

Suppose that $G$ is a minimum counterexample to Theorem~\ref{thm:K9=Extremal}, and put $n := |G|$.
Then $e(G) \ge 6n - 20$, but $G$ is not a $(K_8, K_{2, 2, 2, 2, 2}, 5)$-cockade and $G \not \ge K_9^=$.
We may suppose that $e(G) = 6n - 20$.
It is easy to verify that $n \ge 11$.
\smallskip

\step{\label{K9=Extremal:deltaG7} $\delta(G) \ge 7$.}

Let $x \in V(G)$ such that $d_G(x) = \delta(G)$.
Then $e(G - x) = 6n - 20 - \delta(G) = 6(n - 1) - 14 - \delta(G)$.
If $e(G - x) > 6(n - 1) - 20$, then by the minimality of $G$ and Lemma~\ref{lem:K85CockadeEdges}, $G - x \ge K_9^=$, a contradiction since $G \ge G - x$.
If $e(G - x) = 6(n - 1) - 20$, then $d_G(x) = 6$, and since $G \not\ge K_9^=$, we see by the minimality of $G$ that $G - x$ is a $(K_8, K_{2, 2, 2, 2, 2}, 5)$-cockade.
But then by Lemma~\ref{lem:K85Cockade+v}, $G \ge K_9^=$, a contradiction.
Hence, $e(G - x) < 6(n - 1) - 20$.
It follows that $\delta(G) > 6$.
\proofsquare
\smallskip

\step{\label{K9=Extremal:deltaNx} $\delta(G[N(x)]) \ge 6$ for all $x \in V(G)$.}

Let $xy \in E(G)$ and put $d := |N(x) \cap N(y)|$.
Then $e(G/xy) = 6n - 20 - (d + 1) = 6(n - 1) - 15 - d$.
If $e(G/xy) > 6(n - 1) - 20$, then by Lemma~\ref{lem:K85CockadeEdges} and the minimality of $G$, $G/xy \ge K_9^=$, a contradiction since $G \ge G/xy$.
If $e(G/xy) = 6(n - 1) - 20$, then $d = 5$, and since $G \not\ge K_9^=$, we see by the minimality of $G$ that $G/xy$ is a $(K_8, K_{2, 2, 2, 2, 2}, 5)$-cockade.
But then by Lemma~\ref{lem:K85CockadeContract} and~\stepref{K9=Extremal:deltaG7}, $G \ge K_9^=$, a contradiction.
Hence, $e(G/xy) < 6(n - 1) - 20$.
It follows that $d > 5$.
\proofsquare
\smallskip

We will utilize the following notation throughout the remainder of the proof.
Let $S$ be a minimal separating set in $G$, and let $G_1$ and $G_2$ be two subgraphs of $G$ such that $G = G_1 \cup G_2$ and $G_1 \cap G_2 = G[S]$.
It is an immediate consequence of~\stepref{K9=Extremal:deltaG7} that
\smallskip

\step{\label{K9=Extremal:Gi8} $|G_i| \ge 8$ for $i \in \{1, 2\}$.}
\smallskip

\step{\label{K9=Extremal:GiCockade} Neither $G_1$ nor $G_2$ is a $(K_8, K_{2, 2, 2, 2, 2}, 5)$-cockade.}

Suppose that $G_1$, say, is a $(K_8, K_{2, 2, 2, 2, 2}, 5)$-cockade.
If $|S| \ge 6$, then by contracting any component of $G_2 - S$ to a single vertex, it follows from Lemma~\ref{lem:K85Cockade+v} that $G \ge K_9^=$, a contradiction.
Hence $|S| \le 5$, and it follows from Lemma~\ref{lem:K85Cockade+e} that $G[S]$ is complete.
Then since $e(G_2) = e(G) - e(G_1) + e(G[S])$, from Lemma~\ref{lem:K85CockadeEdges} we have $e(G_2) = (6n - 20) - (6|G_1| - 20) + \binom{|S|}{2} = 6|G_2| - 6|S| + \binom{|S|}{2}$.
If $|S| \le 4$, then $e(G_2) \ge 6|G_2| - 18$, and it follows from~\stepref{K9=Extremal:Gi8}, Lemma~\ref{lem:K85CockadeEdges}, and the minimality of $G$ that $G_2 \ge K_9^=$, a contradiction.
Hence $|S| = 5$, and we have $e(G_2) = 6|G_2| - 20$.
Since $G \not\ge K_9^=$, we see $G_2 \not\ge K_9^=$, and so by~\stepref{K9=Extremal:Gi8} and the minimality of $G$, $G_2$ is a $(K_8, K_{2, 2, 2, 2, 2}, 5)$-cockade.
However, since $G[S]$ is isomorphic to $K_5$, it follows that $G$ is also a $(K_8, K_{2, 2, 2, 2, 2}, 5)$-cockade, a contradiction.
\proofsquare
\smallskip

For $i \in \{1, 2\}$, let $d_i$ be the maximum number of edges that can be added to $G_{3 - i}$ by contracting edges of $G$ with at least one end in $G_i$.
More precisely, let $d_i$ be the largest integer such that $G_i$ contains disjoint sets of vertices $V_1, V_2, \dots, V_{|S|}$ so that $G_i[V_j]$ is connected and $|S \cap V_j| = 1$ for $j \in \{1, \dots, |S|\}$, and such that the graph obtained from $G_i$ by contracting $V_1, V_2, \dots, V_p$ each to a single vertex and deleting $V(G) \setminus \left(\cup_j V_j\right)$ has $e(G[S]) + d_i$ edges.
Let $G_i^*$ be a graph obtained from $G$ by contracting $G_{3 - i}$ onto $S$ so that $G_i^* - S = G_i - S$ and $e(G_i^*) = e(G_i) + d_{3 - i}$ for $i \in \{1, 2\}$.
\smallskip

\step{\label{K9=Extremal:eGi} For $i \in \{1, 2\}$, $e(G_i) \le 6|G_i| - 20 - d_{3 - i}$, with equality only if $G_i^*$ is a $(K_8, K_{2, 2, 2, 2, 2}, 5)$-cockade.}

If, say, $e(G_1) > 6|G_1| - 20 - d_2$, then $e(G_1^*) > 6|G_1^*| - 20$, and so it follows from~\stepref{K9=Extremal:Gi8}, Lemma~\ref{lem:K85CockadeEdges}, and the minimality of $G$ that $G_1^* \ge K_9^=$, a contradiction.
If $e(G_1) = 6|G_1| - 20 - d_2$, then $e(G_1^*) = 6|G_1^*| - 20$, and then since $G \not\ge K_9^=$, it now follows from the minimality of $G$ that $G_1^*$ is a $(K_8, K_{2, 2, 2, 2, 2}, 5)$-cockade.
\proofsquare
\smallskip

\step{\label{K9=Extremal:eG} $e(G) \le 6n - 40 + 6|S| - d_1 - d_2 - e(G[S])$.}

This follows from~\stepref{K9=Extremal:eGi} and the fact that $e(G) = e(G_1) + e(G_2) - e(G[S])$.
\proofsquare
\smallskip

\step{\label{K9=Extremal:5Conn} $G$ is 5-connected.}

Since $e(G) = 6n - 20$, it follows from~\stepref{K9=Extremal:eG} that $6|S| \ge 20 + d_1 + d_2 + e(G[S])$.
Now since $e(G[S]) \ge \frac{1}{2}\delta(G[S])|S|$ and $d_i \ge |S| - \delta(G[S]) - 1$ for $i \in \{1, 2\}$, we see $4|S| \ge 18 + \frac{1}{2}\delta(G[S])(|S| - 4)$.
Therefore $|S| \ge 5$.
\proofsquare
\smallskip

\step{\label{K9=Extremal:d11} If $|S| \in \{5, 6\}$, then for $i \in \{1, 2\}$ there exists $x \in V(G_i) \setminus S$ with $d_G(x) \le 11$.}

Suppose to the contrary that $d_G(x) \ge 12$ for every vertex $x \in V(G_1) \setminus S$, say.
Then $|G_1| \ge 13$ and $|G_1| - |S| \ge 7$.
Then $2e(G_1) \ge 12(|G_1| - |S|) + e(V(G_1) \setminus S, S) + 2e(G[S])$.
By~\stepref{K9=Extremal:eGi}, we also have $2e(G_1) \le 12|G_1| - 40 - 2d_2$.
Hence $e(V(G_1) \setminus S, S) \le 12|S| - 40 - 2(d_2 + e(G[S]))$.
Now there exists $y_1 \in V(G_2) \setminus S$.
Since $|S| \le 6$ and $d_G(y_1) \ge 7$ by~\stepref{K9=Extremal:deltaG7}, $y_1$ has a neighbor $y_2 \in V(G_2) \setminus S$.
By~\stepref{K9=Extremal:deltaNx}, $y_1$ and $y_2$ have at least 6 common neighbors.
Say $y_3, y_4, \dots, y_{|S|}$ are common neighbors of $y_1$ and $y_2$ in $G$.
Note that $N(y_1) \subseteq V(G_2)$.
Now there exist disjoint paths $P_1, \dots, P_{|S|}$ in $G$, each with one end in $\{y_1, \dots, y_{|S|}\}$, the other end in $S$, and all internal vertices in $G_2 - S$.
By contracting each of these paths onto its end in $S$, it follows that $e(G[S]) + d_2 \ge 2(|S| - 2) + 1 = 2|S| - 3$.
Thus $e(V(G_1) \setminus S, S) \le 8|S| - 34$.
Therefore $2e(G_1 - S) \ge 12(|G_1| - |S|) - e(V(G_1) \setminus S, S) \ge 12(|G_1| - |S|) - 8|S| + 34$, and so $e(G_1 - S) \ge 6(|G_1| - |S|) - 7$ since $|S| \in \{5, 6\}$.
Since $e(G_1 - S) \le \binom{|G_1| - |S|}{2}$, this implies $|G_1| - |S| \ge 12$.
From the minimality of $G$ and Lemma~\ref{lem:K85CockadeEdges}, we see $G_1 - S \ge K_9^=$, a contradiction.
\proofsquare
\smallskip

\step{\label{K9=Extremal:6Conn} $G$ is 6-connected.}

Suppose not.
By~\stepref{K9=Extremal:5Conn}, there exists a minimal separating set $S$ of $G$ with $|S| = 5$.
By~\stepref{K9=Extremal:d11}, there exists $x \in V(G_1) \setminus S$ such that $d_G(x) \le 11$.
By~\stepref{K9=Extremal:deltaG7}, $d_G(x) \ge 7$, and by~\stepref{K9=Extremal:deltaNx}, $\delta(G[N(x)]) \ge 6$.
If $G[N(x)]$ contains a subgraph $H$ isomorphic to $K_5$, then let $P_1, \dots, P_5$ be disjoint paths in $G$ with one end in $V(H)$ and one end in $S$.
Otherwise, let $P_1, \dots, P_5$ be disjoint paths with one end in $N(x)$, the other end in $S$, and no internal vertices in $N(x)$.
Say $V(P_i) \cap N(x) = v_i$ for $i \in \{1, \dots, 5\}$.
Then by Lemma~\ref{lem:11VK5} the end of some $P_i$ in $N(x)$, say $v_1$, is such that $G[N(x)] - v_1$ has a $K_4$ minor rooted at $\{v_2, v_3, v_4, v_5\}$.
By additionally contracting the edge $xv_1$ if necessary, we may thus in both cases assume that the ends of $P_1, \dots, P_5$ in $N(x)$ belong to a $K_5$ subgraph of some minor of $G[N[x]]$.
By now contracting each path $P_i$ onto its end in $S$, it follows that $d_1 + e(G[S]) = 10$.
Similarly, $d_2 + e(G[S]) = 10$, and so $d_1 = d_2$.
Note that by~\stepref{K9=Extremal:GiCockade} and~\stepref{K9=Extremal:eGi}, for $i \in \{1, 2\}$ we must have either $d_{3 - i} \ge 1$ or $e(G_i) \le 6|G_i| - 21 - d_{3 - i}$.
Then in either case, $6n - 20 = e(G) = e(G_1) + e(G_2) - e(G[S]) \le 6n - 10 - \max\{1, d_1\} - (d_2 + e(G[S])) \le 6n - 21$, a contradiction.
\proofsquare
\smallskip

\step{\label{K9=Extremal:CliqueSep} There is no minimal separating set $S$ of $G$ such that $G[S]$ is complete.}

Suppose such a separating set $S$ exists.
If $|S| \ge 7$, then by contracting any two components of $G - S$ to a single vertex each, we see $G \ge K_9^-$, a contradiction.
Thus, by~\stepref{K9=Extremal:6Conn}, $|S| = 6$.
Since $G[S]$ is complete, we have $d_1 = d_2 = 0$.
For $i \in \{1, 2\}$, since $G_i \not \ge K_9^=$ and $G_i$ is not a $(K_8, K_{2, 2, 2, 2, 2}, 5)$-cockade by~\stepref{K9=Extremal:GiCockade}, we have $e(G_i) \le 6|G_i| - 21$ by~\stepref{K9=Extremal:eGi}.
But then $6n - 20 = e(G) = e(G_1) + e(G_2) - e(G[S]) \le (6|G_1| - 21) + (6|G_2| - 21) - 15 = 6n - 21$, a contradiction.
\proofsquare
\smallskip

\step{\label{K9=Extremal:deltaG8} $\delta(G) \ge 8$.}

Suppose there exists $x \in V(G)$ such that $d_G(x) \le 7$.
By~\stepref{K9=Extremal:deltaG7}, $d_G(x) = 7$.
By~\stepref{K9=Extremal:deltaNx}, $\delta(G[N(x)]) = 6$, and so $G[N(x)]$ is isomorphic to $K_7$.
But $N(x)$ is a separating set of $G$, and so contains some minimal separating set of $G$, contrary to~\stepref{K9=Extremal:CliqueSep}.
\proofsquare
\smallskip

\step{\label{K9=Extremal:CliqueMinusSep} There is no minimal separating set $S$ of $G$ such that $G[S \setminus \{x\}]$ is complete for some $x \in S$.}

Suppose there exists $x \in S$ such that $G[S \setminus \{x\}]$ is complete.
Note that $|S| \ge 6$ by~\stepref{K9=Extremal:6Conn}.
If $|S| \ge 8$, then since $G[S \setminus \{x\}]$ contains $K_7$ as a subgraph, we can contract any two components of $G - S$ each to a single vertex to obtain a $K_9^-$ minor, a contradiction.
So $|S| \in \{6, 7\}$.
By contracting any component of $G_i - S$ onto $x$, we see that $d_i + e(G[S]) = \binom{|S|}{2}$ for $i \in \{1, 2\}$, and in particular $d_1 = d_2 = |S| - \delta(G[S]) - 1 \ge 1$.
Suppose, say, $e(G_1) \le 6|G_1| - 21 - d_2$.
Then $6n - 20 = e(G) = e(G_1) + e(G_2) - e(G[S]) \le (6|G_1| - 21 - d_2) + (6|G_2| - 20 - d_1) - e(G[S]) \le 6n - 21$, a contradiction.

Hence by~\stepref{K9=Extremal:eGi}, we have $e(G_i) = 6|G_i| - 20 - d_{3 - 1}$, and $G_i^*$ is a $(K_8, K_{2, 2, 2, 2, 2}, 5)$-cockade for $i \in \{1, 2\}$.
Now $6n - 20 = e(G) = e(G_1) + e(G_2) - e(G[S]) = (6|G_1| - 20 - d_2) + (6|G_2| - 20 - d_1) - e(G[S]) = 6n - 19 - d_1$.
Thus we must have $d_1 = d_2 = 1$, and so $G[S]$ is isomorphic to $K_{|S|}^-$.
Since $G_i^*[S]$ contains a $K_6$ subgraph, there must exist some $K_8$ subgraph $H_i$ of $G_i^*$ such that $S \subseteq V(H_i)$ for $i \in \{1, 2\}$.
Then $G[V(H_i)]$ is isomorphic to $K_8^-$.
Thus $G[V(H_1 \cup H_2)]$ is not a $(K_8, K_{2, 2, 2, 2, 2}, 5)$-cockade, but is a graph on $16 - |S|$ vertices with $6(16 - |S|) - 20$ edges.
Since $n \ge 11 > 16 - |S|$, $G[V(H_1 \cup H_2)]$ contains a $K_9^=$ minor by the minimality of $G$, a contradiction.
\proofsquare
\smallskip

\step{\label{K9=Extremal:deltaG9} $\delta(G) \ge 9$.}

Suppose there exists $x \in V(G)$ such that $d_G(x) \le 8$.
By~\stepref{K9=Extremal:deltaG8}, $d_G(x) = 8$.
Since $\delta(G[N(x)]) \ge 6$ by~\stepref{K9=Extremal:deltaNx}, it follows that $G[N(x)]$ is isomorphic to $K_8 - M$, where $M$ is a matching of $N(x)$.
If $|M| \le 2$, then $G[N[x]]$ contains $K_9^=$ as a subgraph, a contradiction.
Suppose $|M| = 3$.
Let $C$ be a component of $G - N[x]$.
Since $G$ is 6-connected by~\stepref{K9=Extremal:6Conn}, $N(C)$ contains both ends of some missing edge $e$ of $G[N(x)]$.
By contracting $C$ onto one end of $e$, we see that $G \ge K_9^=$, a contradiction.
Hence $|M| = 4$.
That is, $G[N(x)]$ is isomorphic to $K_{2, 2, 2, 2}$.
Say $N(x) = \{v_1, w_1, \dots, v_4, w_4\}$, where $v_i w_i \notin E(G)$ for $i \in \{1, \dots, 4\}$.
\medskip

\sstep{\label{K9=Extremal:deltaG9:OneComp} $G - N[x]$ is connected.}

If not, let $C_1, C_2$ be distinct components of $G - N[x]$.
Since $G$ is 6-connected by~\stepref{K9=Extremal:6Conn}, $N(C_1)$ and $N(C_2)$ each contain at least two nonadjacent pairs of vertices of $N(x)$.
Thus it is possible to pick distinct pairs from each of $N(C_1)$ and $N(C_2)$, say $v_i, w_i \in N(C_i)$ for $i \in \{1, 2\}$.
By contracting $C_i$ onto $v_i$ for $i \in \{1, 2\}$ we see $G \ge K_9^=$, a contradiction.
\proofsquare
\smallskip

\sstep{\label{K9=Extremal:deltaG9:NonadjacentPair} No vertex of $G - N[x]$ is adjacent to a pair of nonadjacent vertices of $N(x)$.}

Suppose to the contrary that there exists $z \in V(G) \setminus N[x]$ such that $z$ is adjacent to, say, $v_1$ and $w_1$.
Let $C$ be a component of $G - N[x] - z$.
If $N(C)$ contains some pair of nonadjacent vertices of $N(x)$ distinct from $v_1, w_1$, say $v_2, w_2 \in N(C)$, then by contracting $C$ onto $v_2$ and contracting the edge $v_1z$, we see $G \ge K_9^=$, a contradiction.
Hence $N(C)$ can contain no nonadjacent pair of vertices of $N(x)$ other than $v_1, w_1$.
It follows that $|N(C) \cap N(x)| \le 5$.
Since $G$ is 6-connected by~\stepref{K9=Extremal:6Conn}, we conclude $|N(C) \cap N(x)| = 5$, and $\{v_1, w_1, z\} \subseteq N(C)$.
Without loss of generality, $N(C) = \{v_1, v_2, v_3, v_4, z, w_1\}$.
Since $N(C)$ is a minimal separating set of $G$, we may put $S = N(C)$.
Note that $G[\{v_1, w_1, v_2, v_3, v_4\}]$ is isomorphic to $K_5^-$, so $e(G[S]) \ge 11$.
By~\stepref{K9=Extremal:d11}, there exists $y \in V(G_1) \setminus S$ such that $d_G(y) \le 11$.
Since $G$ is $6$-connected by~\stepref{K9=Extremal:6Conn}, there exist disjoint paths $P_1, \dots, P_6$ with one end in $N(y)$, the other end in $S$, and no internal vertices in $N(y) \cup S$.
Say $u_i$ is the end of $P_i$ in $N(y)$ for $i \in \{1, \dots, 6\}$, $v_i$ is the end of $P_i$ in $S$ for $i \in \{1, \dots, 4\}$, and $P_5, P_6$ have ends $z, w_1$ in $S$, respectively.
If $u_1 u_6 \in E(G)$, then by contracting each path $P_i$ onto its end in $S$ and additionally contracting the edge $y u_5$, we see $e(G[S]) + d_2 = 15$.
Otherwise, by Lemma~\ref{lem:EdgeorK14} there exists some component $C'$ of $G[N(y)] - \{u_1, \dots, u_6\}$ such that either $\{u_1, u_6\} \subseteq N(C')$ or $\{u_2, \dots, u_5\} \subseteq N(C')$.
By contracting $C'$ onto $u_1$ and contracting $y u_5$ in the former case, or contracting $C'$ onto $u_5$ and contracting $y u_1$ in the latter case, and then contracting each path $P_i$ onto its end in $S$ we again see $e(G[S]) + d_2 = 15$.
By symmetry, $e(G[S]) + d_1 = 15$ as well.
By~\stepref{K9=Extremal:CliqueMinusSep}, $d_i \ge 2$ for $i \in \{1, 2\}$.
Thus by~\stepref{K9=Extremal:eG}, we see $6n - 20 = e(G) \le 6n - 40 + 36 - (e(G[S]) + d_2) - d_1 \le 6n - 21$, a contradiction.
\proofsquare
\smallskip

\sstep{\label{K9=Extremal:deltaG9:ThreeNeighbors} Every vertex of $N(x)$ has at least three neighbors in $G - N[x]$.}

Since $\delta(G) \ge 8$ by~\stepref{K9=Extremal:deltaG8}, there exists $z \in V(G) \setminus N[x]$ such that $z v_1 \in E(G)$, say.
By~\stepref{K9=Extremal:deltaNx}, $z$ and $v_1$ have at least six common neighbors.
By~\sstepref{K9=Extremal:deltaG9:NonadjacentPair}, at most three of these common neighbors can belong to $N(x)$, and~\sstepref{K9=Extremal:deltaG9:ThreeNeighbors} follows.
\proofsquare
\smallskip

\sstep{\label{K9=Extremal:deltaG9:2Conn} $G - N[x]$ is 2-connected.}

By~\sstepref{K9=Extremal:deltaG9:OneComp}, suppose to the contrary that $z$ is a cut-vertex of $G - N[x]$.
If $C_1$ and $C_2$ are components of $G - N[x] - z$ with, say, $v_i, w_i \in N(C_i)$ for $i \in \{1, 2\}$, then by contracting each $C_i$ onto $v_i$ we see $G \ge K_9^=$, a contradiction.
Since $G$ is $6$-connected by~\stepref{K9=Extremal:6Conn}, it follows that every component $C$ of $G - N[x] - z$ satisfies $|N(C) \cap N(x)| = 5$, with $v_1$ and $w_1$, say, in $N(C)$, and at most one of $v_i$ or $w_i$ in $N(C)$ for $i \in \{2, 3, 4\}$.
Without loss of generality, assume $C_1$ is a component of $G - N[x] - z$ with $\{v_2, v_3, v_4\} \subseteq N(C)$.
Consider the three edges $e_1 = v_2 w_3$, $e_2 = v_3 w_2$, and $e_3 = w_2 w_3$.
Then for $i \in \{1, 2, 3\}$, by~\stepref{K9=Extremal:deltaNx} the ends of $e_i$ have at least six common neighbors, only five of which belong to $N[x]$, and so there exists $y_i \in V(G) \setminus N[x]$ adjacent to both ends of $e_i$.
By~\sstepref{K9=Extremal:deltaG9:NonadjacentPair}, $y_1, y_2, y_3$ are distinct, and in particular at most one $y_i = z$.
As $w_2, w_3 \notin N(C_1)$, $y_i \notin V(C_1)$ for $i \in \{1, 2, 3\}$.
Furthermore, if two of $y_1, y_2, y_3$ belong to the same component $C_2$ of $G - N[x] - z$, then $C_2$ is a component with $v_2, w_2 \in N(C_2)$, say, contrary to the above.
Thus there exist at least two components $C_2$ and $C_3$ of $G - N[x] - z$ distinct from $C_1$.
Without loss of generality, we may assume $w_2 \in N(C_2)$.
Now by contracting $C_1$ and $C_2$ onto $z$, contracting $C_3$ onto $v_1$, and contracting the edge $v_2 z$, we see $G \ge K_9^=$, a contradiction.
\proofsquare
\smallskip

We now consider the graph $H = G - \{x, v_3, w_3, v_4, w_4\}$.
We claim that $H$ is 4-connected.
Suppose $Q$ is a minimum separating set of at most three vertices in $H$.
By~\sstepref{K9=Extremal:deltaG9:2Conn}, we have $|Q| \ge 2$ and $|Q \cap N(x)| \le 1$.
If $|Q \cap N(x)| = 1$, then by symmetry we may assume $w_2 \in Q$.
Since $H[\{v_1, w_1, v_2\}]$ is connected, $v_1$, $w_1$, and $v_2$ all belong to the same component $C$ of $H - Q$.
If $w_2 \notin Q$, then $w_2$ also belongs to $C$, and in this case we assume that $Q$ and $w_2$ are chosen so that $|Q \cap N(w_2)|$ is maximal.
We next claim that there exist $v_1'$ and $w_1'$ in $V(G) \setminus (N[x] \cup Q)$ adjacent to $v_1$ and $w_1$, respectively.
If not, then by~\sstepref{K9=Extremal:deltaG9:NonadjacentPair} and~\sstepref{K9=Extremal:deltaG9:ThreeNeighbors}, it must be the case that $w_1$, say, has exactly three neighbors $z_1, z_2, z_3$ in $G - N[x]$, and $Q = \{z_1, z_2, z_3\}$.
Now $w_2 \notin Q$, so by our choice of $Q$ and $w_2$, it follows that $w_2$ is complete to $Q$.
Since $v_2 w_1 \in E(G)$, $v_2$ and $w_1$ have at least one common neighbor in $G - N[x]$ by~\stepref{K9=Extremal:deltaNx}.
This common neighbor must be one of $z_1, z_2, z_3$, say $z_1$, but then $z_1$ is adjacent to $v_2$ and $w_2$, contradicting~\sstepref{K9=Extremal:deltaG9:NonadjacentPair}.
Thus the claim is proved, and there exist $v_1', w_1' \in V(G) \setminus (N[x] \cup Q)$ such that $v_1 v_1', w_1 w_1' \in E(G)$.
Now we have $v_1', w_1' \in V(C)$.
By~\sstepref{K9=Extremal:deltaG9:2Conn}, there exist two internally disjoint $v_1', w_1'$-paths in $G - N[x]$.
Since $|Q| \le 3$, at least one of these paths must be contained entirely within $G[V(C) \cup Q]$.
Note that since $G \not\ge K_9^=$, there must then be no $v_i, w_i$-path in $G[V(C') \cup \{v_i, w_i\}]$ for $i \in \{3, 4\}$, where $C'$ is any component of $H - Q$ distinct from $C$.
Hence at most one of $v_i, w_i$ has a neighbor in $C'$ for $i \in \{3, 4\}$.
It follows that $C'$ is separated from $x$ by $Q$ and at most two vertices of $N(x)$.
But since $|Q| \le 3$, this contradicts that $G$ is 6-connected by~\stepref{K9=Extremal:6Conn}.
This proves the claim that $H$ is 4-connected.

If there exists a $K_4$ minor of $H$ rooted at $\{v_1, w_1, v_2, w_2\}$, then $G \ge K_9^=$, a contradiction.
Thus $e(H) \le 3|H| - 7 = 3(n - 5) - 7$ by Theorem~\ref{thm:RootedK4}.
For $i \in \{3, 4\}$, $v_i$ and $w_i$ have no common neighbor in $G - N[x]$ by~\sstepref{K9=Extremal:deltaG9:NonadjacentPair}, so they together have at most $|G| - |N[x]| = n - 9$ neighbors in $G - N[x]$.
Furthermore, the vertices $v_3, w_3, v_4, w_4$ are together incident with 20 edges of $G[N(x)]$.
Therefore $6n - 20 = e(G) \le d_G(x) + 20 + 2(n - 9) + e(H) \le 8 + 20 + 2(n - 9) + 3(n - 5) - 7 = 5n - 12$.
It follows that $n \le 8$, a contradiction which completes the proof of~\stepref{K9=Extremal:deltaG9}.
\proofsquare
\smallskip

\step{\label{K9=Extremal:NbrComp} Suppose $x \in V(G)$ with $d_G(x) \in \{9, 10, 11\}$, and let $M \subseteq N(x)$ be the vertices of $N(x)$ which are not complete to all other vertices of $N(x)$.
Then there is no component $C$ of $G - N[x]$ such that $N(C') \cap M \subseteq N(C)$ for all components $C'$ of $G - N[x]$.}

Suppose to the contrary that such a component $C$ exists.
Among all vertices $x$ with $d_G(x) \in \{9, 10, 11\}$ for which such a component $C$ exists, choose $x$ to be of minimum degree.
Note that $N(C') \cap M \ne \emptyset$ for all components $C'$ of $G - N[x]$ by~\stepref{K9=Extremal:CliqueSep}.
Suppose for a contradiction that $M \setminus N(C) \ne \emptyset$, and choose $y \in M \setminus N(C)$ to be of minimum degree among all vertices in $M \setminus N(C)$.
Then $d_G(y) < d_G(x)$ since $y$ has no neighbor outside $N[x]$ by the existence of $C$.
Now let $K$ be the component of $G - N[y]$ containing $C$.
We claim that $N(x) \setminus N[y] \not\subseteq V(K)$.
So suppose instead that $N(x) \setminus N[y] \subseteq V(K)$.
Let $C'$ be any component of $G - N[x]$ distinct from $C$.
Note that such a component $C'$ exists since otherwise $K$ is the only component of $G -N[y]$, contrary to our choice of $x$ and $C$.
Let $K'$ be the component of $G - N[y]$ containing $C'$.
We may assume $C'$ is chosen such that $K' \ne K$, since otherwise $K$ is the only component of $G - N[y]$, again a contradiction.
Then $V(K') \cap (N(x) \setminus N[y]) = \emptyset$, since $N(x) \setminus N[y] \subseteq V(K)$.
Hence $N(K') = N(C') \subseteq N(y)$.
Thus we have that $N(C') \cap M \subseteq N(C) \cap N(y)$.
Therefore $N(K') \cap M_y \subseteq N(K)$, where $M_y$ is the set of vertices of $N(y)$ not complete to all other vertices of $N(y)$.
Noticing that $M_y \subseteq M$ and that the component $K'$ was essentially arbitrary (every component $K'$ of $G - N[y]$ corresponds with some component $C'$ of $G - N[x]$), we see that the existence of $y$ and $K$ contradicts the choice of $x$ and $C$.
Therefore $N(x) \setminus N[y] \not\subseteq V(K)$, as claimed.

Hence there exists some component $H$ of $G[N(x) \setminus N[y]]$ with $V(H) \cap N(C) = \emptyset$.
We must have $d_G(z) \ge d_G(y)$ for all $z \in V(H)$ by the choice of $y$.
If $|H| = 1$, then it follows that $d_G(z) = d_G(y)$ and $N(z) = N(y)$.
But then $H$ is a component of $G - N[y]$ contradicting the choice of $x$ and $C$.
Thus $|H| \ge 2$.
On the other hand, $|H| \le d_G(x) - d_G(y) \le 11 - 9 = 2$, and so $|H| = 2$, $d_G(x) = 11$, and $d_G(y) = 9$.
By~\stepref{K9=Extremal:deltaNx} applied to $y$, we see that $G[N(y) \cap N(x)]$ has minimum degree at least 5.
Then the edges of $G[N(x)]$ are the edges of $N(x) \cap N(y)$, edges incident with $y$, and edges incident with $V(H)$.
Therefore, $e(G[N(x)]) \ge \frac{5}{2}(d_G(y) - 1) + (d_G(y) - 1) + \left(2(d_G(y) - 1) - 1\right) = 43 > 5d_G(x) - 14$.
By Theorem~\ref{thm:K8=Extremal}, we see that $G[N(x)] \ge K_8^=$, and therefore $G[N[x]] \ge K_9^=$, a contradiction.
This proves that $M \setminus N(C) = \emptyset$, that is $M \subseteq N(C)$.

If $G[N(x)] \ge K_7^= \cup K_1$, then let $y \in N(x)$ such that $G[N(x) \setminus \{y\}] \ge K_7^=$.
If $y \notin N(C)$, then $y$ is complete to $N(x) \setminus \{y\}$, and so $G[N(x)] \ge K_8^=$ and $G[N[x]] \ge K_9^=$, a contradiction.
Hence $y \in N(C)$ and every nonneighbor of $y$ in $N(x)$ also belongs to $N(C)$ since $M \subseteq N(C)$.
Now by contracting $C$ onto $y$, we again find, along with $x$, a $K_9^=$ minor in $G$, a contradiction.
Hence by Lemma~\ref{lem:K7=UK1Computer}, $G[N(x)]$ is isomorphic to one of the five graphs $\overline{C_5} \vee \overline{C_4}$, $\overline{C_9}$, $K_{3, 3, 3}$, $\overline{C_6} \vee \overline{K_3}$, or $\overline{P}$, where $\overline{P}$ is the complement of the Petersen graph.
Suppose $G[N(x)]$ is isomorphic to one of $\overline{C_5} \vee \overline{C_4}$ or $\overline{C_9}$.
In both cases, $G[N(x)] \ge K_7^-$ and $N(x) = M \subseteq N(C)$, so by contracting $C$ to a single vertex we obtain, along with $x$, a $K_9^=$ minor of $G$, a contradiction.
Thus we may suppose $G[N(x)]$ is isomorphic to one of $K_{3, 3, 3}$, $\overline{C_6} \vee \overline{K_3}$, or $\overline{P}$.
Note that by Lemma~\ref{lem:K7=UK1Computer}, these three graphs are edge-maximal subject to not having a $K_7^= \cup K_1$ minor.
We first show the following.
\medskip

\sstep{\label{K9=Extremal:NbrComp:OneComp} $G - N[x]$ is connected.}

Suppose $C'$ is a component of $G - N[x]$ distinct from $C$.
By~\stepref{K9=Extremal:CliqueSep}, $N(C')$ contains both ends of some missing edge $e$ of $G[N(x)]$, and we contract $C'$ onto one end of $e$.
Then $G[N(x)] + e \ge K_7^= \cup K_1$.
Now let $y \in N(x)$ such that $G[N(x) \setminus \{y\}] + e \ge K_7^=$.
By contracting $C$ onto $y$, we see $G \ge K_9^=$.
\proofsquare
\smallskip

\sstep{\label{K9=Extremal:NbrComp:K333} $G[N(x)]$ is not isomorphic to either $K_{3, 3, 3}$ or $\overline{C_6} \vee \overline{K_3}$.}

Say $N(x) = \{v_1, \dots, v_9\}$ where $\{v_1, v_2, v_3\}$ is an independent set and either $\{v_4, v_5, v_6\}$ and $\{v_7, v_8, v_9\}$ are independent sets if $G[N(x)]$ is isomorphic to $K_{3, 3, 3}$ or $\{v_4, \dots, v_9\}$ are the vertices of a $C_6$ in $\overline{G[N(x)]}$ written in cyclic order if $G[N(x)]$ is isomorphic to $\overline{C_6} \vee \overline{K_3}$.
We claim that $G - N[x]$ is 3-connected.
Suppose to the contrary that $Q$ is a minimum cut set of $G - N[x]$ with $|Q| \le 2$.
By~\sstepref{K9=Extremal:NbrComp:OneComp}, $|Q| \ge 1$.
Let $C_1$ be any component of $G - (N[x] \cup Q)$, and let $C_2 := G - (N[x] \cup V(C_1))$.
Then $C_2$ is connected.
Suppose that some vertex $v_i \in N(x)$ has no neighbor in $V(C_2)$.
By~\stepref{K9=Extremal:deltaNx}, any neighbor of $v_i$ in $N(x)$ must have at least 2 neighbors in $V(C_1)$, that is, $\{v_j \in N(x) : v_i v_j \in E(G) \} \subseteq N(C_1)$.
In particular, there exist disjoint sets $T_1, T_2 \subseteq N(C_1)$ such that $|T_k| = 3$ and $G[T_k]$ contains at least two missing edges of $G[N(x)]$ for $k \in \{1, 2\}$.
Now since $C_2$ contains at least one component of $G - (N[x] \cup Q)$, it follows by~\stepref{K9=Extremal:6Conn} that $|N(C_2) \cap N(x)| \ge 4$.
If $G[N(C_2)]$ contains some missing edge $e$ of $G[N(x)]$, then some $T_k$, say $T_1$, does not contain both ends of $e$.
We contract $C_2$ onto one end of $e$, and we contract $C_1$ onto a vertex in $T_1$ incident to two missing edges of $G[N(x)]$, and from Lemma~\ref{lem:K333} we see, along with $x$, that $G \ge K_9^=$, a contradiction.
If $G[N(C_2)]$ does not contain a missing edge of $G[N(x)]$, then it must be the case that $G[N(x)]$ is isomorphic to $\overline{C_6} \vee \overline{K_3}$, $|Q| = 2$, and $N(C_2) \cap N(x) = \{v_3, v_4, v_6, v_8\}$, say.
Then there exist distinct $w_1, w_2 \in V(C_1)$ such that $w_1 v_1, w_2 v_2 \in E(G)$, and disjoint paths $P_1, P_2$ in $G - N[x]$ where $P_i$ has one end $w_i$ and one end in $Q$.
Furthermore, since $C_1$ is connected there exists a path $P_3$ in $C_1$ with one end in $P_1$, the other end in $P_2$, and no internal vertices in $P_1 \cup P_2$.
Now by contracting $C_2$ onto $v_3$, contracting the paths $P_i \cup w_i v_i$ onto $v_i$ for $i \in \{1, 2\}$, and contracting $P_3$ to a single edge, we see $G \ge K_9^=$ by Lemma~\ref{lem:K333}, a contradiction.
Therefore we may assume that every vertex of $N(x)$ has some neighbor in $V(C_2)$.
But now $|N(C_1) \cap N(x)| \ge 4$, and the same argument as above with the roles of $C_1$ and $C_2$ switched will give a contradiction.
This establishes that $G - N[x]$ is 3-connected.

Note that since $\delta(G) \ge 9$ by~\stepref{K9=Extremal:deltaG9}, every vertex of $N(x)$ has at least two neighbors in $G - N[x]$.
If some vertex $z \in V(G) \setminus N[x]$ is adjacent to both ends of a missing edge $e$ of $G[N(x)]$, then by contracting $G - N[x] - z$ onto a vertex of $N(x)$ incident to two missing edges of $G[N(x)]$ distinct from $e$, and then contracting $z$ onto one end of $e$, we get a $K_9^=$ minor by Lemma~\ref{lem:K333}, a contradiction.
Hence no vertex of $V(G) \setminus N[x]$ is adjacent to both ends of a missing edge of $G[N(x)]$.
Thus we may select distinct $z_1, z_2, z_3 \in V(G) \setminus N[x]$ such that $v_i z_i \in E(G)$ for $i \in \{1, 2, 3\}$.
Since $G - N[x]$ is 3-connected, by an application of Menger's Theorem there exists a cycle in $G - N[x]$ containing all of $z_1, z_2, z_3$.
By contracting each of the three subpaths of the cycle between the $z_i$ to a single edge, and then contracting each of the three edges $v_i z_i$, we again have $G \ge K_9^=$ by Lemma~\ref{lem:K333}.
This contradiction proves~\sstepref{K9=Extremal:NbrComp:K333}.
\proofsquare
\smallskip

Hence we may assume $G[N(x)]$ is isomorphic to $\overline{P}$, and we label the vertices of $\overline{G[N(x)]}$ as in Figure~\ref{fig:Petersen}.
\smallskip

\sstep{\label{K9=Extremal:NbrComp:PNbr} No vertex of $G - N[x]$ is adjacent to both ends of a missing edge of $G[N(x)]$.}

Suppose there exists $z \in V(G) \setminus N[x]$ adjacent to both ends of some missing edge $e$ of $G[N(x)]$, say $zv_0, zv_1 \in E(G)$.
Note that every vertex of $N(x)$ has at least two neighbors in $G - N[x]$ since $\delta(G) \ge 9$ by~\stepref{K9=Extremal:deltaG9}.
Thus if $G - N[x] - z$ is connected, then by contracting $z$ onto $v_0$ and contracting $G - N[x] - z$ onto $v_3$, say, we add four edges to $G[N(x)]$, and it follows from Lemma~\ref{lem:Petersen} that $G \ge K_9^=$, a contradiction.
So suppose $G - N[x] - z$ is disconnected.
Since $G$ is $6$-connected by~\stepref{K9=Extremal:6Conn}, we have $|N(C) \cap N(x)| \ge 5$ for all components $C$ of $G - N[x] - z$, and so $N(C)$ contains both ends of at least two missing edges of $G[N(x)]$.
Let $C_1$ be a component of $G - N[x] - z$, and suppose $e_1 \ne e$ is a missing edge of $G[N(x)]$ with both ends in $N(C_1)$.
We consider three cases.

First, suppose $e$ and $e_1$ belong to a $6$-cycle, but not a $5$-cycle of $\overline{G[N(x)]}$, say $e_1 = v_7 v_9$.
If $v_2 \in N(C_1)$, then by contracting $C_1$ onto $v_7$ and contracting $zv_0$, we see $G \ge K_9^=$ by Lemma~\ref{lem:Petersen}, a contradiction.
Hence $v_2 \notin N(C_1)$, so there exists some component $C_2$ of $G - N[x] - z$ such that $v_2 \in N(C_2)$.
Now since $|N(C_2) \cap N(x)| \ge 5$, there must exist some missing edge $e_2$ of $G[N(x)]$ with both ends in $N(C_2)$ such that $e_2$ is distinct from $e$ and $e_1$.
By contracting $C_2$ onto one end of $e_2$, $C_1$ onto $v_7$, and $z$ onto $v_0$, we again see $G \ge K_9^=$ by Lemma~\ref{lem:Petersen}, a contradiction.

Next, suppose $e$ and $e_1$ share a common end, say $e_1 = v_1 v_2$.
If $v_3 \in N(C_1)$, then by contracting $C_1$ onto $v_2$ and contracting $zv_0$, we see $G \ge K_9^=$ by Lemma~\ref{lem:Petersen}, a contradiction.
Thus there exists a component $C_2 \ne C_1$ of $G - N[x] - z$ such that $v_3 \in N(C_2)$.
From Lemma~\ref{lem:Petersen}, we may assume that $N(C_2)$ does not contain both ends of any missing edges of $G[N(x)]$ other than $v_3 v_4$, $v_5 v_7$, or $v_1 v_6$.
This requires $N(C_2) \cap N(x) \subseteq \{v_1, v_3, v_4, v_5, v_6, v_7\}$.
By relabelling if necessary, we may assume $v_3, v_4 \in N(C_2)$.
By Lemma~\ref{lem:Petersen}, we may assume that $N(C_1)$ does not contain both ends of any missing edges of $G[N(x)]$ other than $e_1$ or $v_0 v_1$.
This requires $N(C_1) = \{v_0, v_1, v_2, v_8, v_9, z\}$.
Let $S = N(C_1)$.
Then $S$ is a minimum separating set of $G$.
Let $G_1$ and $G_2$ be as defined prior to~\stepref{K9=Extremal:Gi8}, where we may assume $V(G_1) \setminus S = V(C_1)$.
By contracting $V(C_2) \cup \{v_3, v_4\}$ onto $z$ and contracting $x v_1$, we see $e(G[S]) + d_2 = 15$.
By contracting $C_1$ onto $v_1$, we see $d_1 \ge 2$.
Thus by~\stepref{K9=Extremal:eG}, we have $6n - 20 = e(G) \le 6n - 40 + 6 \cdot 6 - 2 - 15 = 6n - 21$, a contradiction.

Lastly, suppose $e$ and $e_1$ are disjoint and belong to a $5$-cycle of $G[N(x)]$, say $e_1 = v_2 v_3$.
If $v_7 \in N(C_1)$, then by contracting $C_1$ onto $v_2$ and contracting $z v_0$, $G \ge K_9^=$ by Lemma~\ref{lem:Petersen}, a contradiction.
Thus there exists a component $C_2 \ne C_1$ of $G - N[x] - z$ such that $v_7 \in N(C_2)$.
If any of $v_2$, $v_5$, or $v_9$ belongs to $N(C_2)$, it is now possible to contract $C_1$, $C_2$, and $z$ onto $N(x)$ such that $G \ge K_9^=$ by Lemma~\ref{lem:Petersen}, a contradiction.
So $N(C_2)$ contains at least four vertices of $\{v_0, v_1, v_3, v_4, v_6, v_8\}$.
From Lemma~\ref{lem:Petersen}, we may assume $N(C_2) \cap N(x) = \{v_0, v_1, v_4, v_7, v_8\}$.
But now considering $C_2$ and the missing edge $v_0 v_4$ of $G[N(x)]$ puts us in the previous case.
\proofsquare
\smallskip

\sstep{\label{K9=Extremal:NbrComp:P2Conn} $G - N[x]$ is $2$-connected.}

Suppose that $z$ is a cut-vertex of $G - N[x]$.
We will show that $z$ must be adjacent to both ends of some missing edge of $G[N(x)]$, contrary to~\sstepref{K9=Extremal:NbrComp:PNbr}.
Let $C_1$ be a component of $G - N[x] - z$, chosen such that $|N(C_1) \cap N(x)|$ is minimum among all components of $G - N[x] - z$.
If $|N(C_1) \cap N(x)| \ge 7$, then it follows from Lemma~\ref{lem:Petersen7} that for any component $C_2 \ne C_1$ of $G - N[x] - z$, we may contract $C_1$ and $C_2$ onto $N(x)$ so that $G \ge K_9^=$, a contradiction.
Thus $|N(C_1) \cap N(x)| \le 6$.
Since $G$ is $6$-connected by~\stepref{K9=Extremal:6Conn}, we have $|N(C_1) \cap N(x)| \ge 5$.
Then $N(C_1)$ contains both ends of some missing edge of $G[N(x)]$, say $v_0, v_1 \in N(C_1)$.
Let $e$ be any missing edge of $G[N(x)]$ such that $v_0 v_1$ and $e$ are disjoint and belong to the same $5$-cycle of $\overline{G[N(x)]}$, and suppose neither end of $e$ belongs to $N(C_1)$.
We may assume $e = v_2 v_3$.
Since $\delta(G) \ge 9$ by~\stepref{K9=Extremal:deltaG9}, $v_2$ and $v_3$ each have at least two neighbors in $G - N[x] - V(C_1)$.
Let $w_i$ be a neighbor of $v_i$ in $G - N[x] - V(C_1)$ for $i \in \{2, 3\}$.
We may assume $w_3 \ne z$.
Furthermore, $w_2 \ne w_3$ by~\sstepref{K9=Extremal:NbrComp:PNbr}.
Since $G - N[x] - V(C_1)$ is connected, there exists a path $P_1$ in $G - N[x] - V(C_1)$ with ends $w_2, w_3$, and a path $P_2$ with one end $z$, the other end in $V(P_1)$, and no internal vertices in $V(P_1)$.
Possibly, $P_2$ consists of only the vertex $z$.
Now by contracting $P_1 \cup P_2$ to a single vertex, contracting the edge $v_2 w_2$, and contracting $C_1$ onto $v_1$, we add the edges $v_0 v_1, v_1 v_2, v_2 v_3$ to $G[N(x)]$, and it follows from Lemma~\ref{lem:Petersen} that $G \ge K_9^=$, a contradiction.

Thus at least one end of every missing edge of $G[N(x)]$ disjoint from $v_0 v_1$ and belonging to a $5$-cycle of $\overline{G[N(x)]}$ with $v_0 v_1$ must belong to $N(C_1)$.
There are eight such missing edges of $G[N(x)]$ which give the $8$-cycle $v_2 v_3 v_4 v_9 v_6 v_8 v_5 v_7$ of $\overline{G[N(x)]}$.
Since $|N(C_1) \cap N(x)| \le 6$, either $N(C_1) = \{z, v_0, v_1, v_2, v_4, v_5, v_6\}$ or $N(C_1) = \{z, v_0, v_1, v_3, v_7, v_8, v_9\}$.
In either case, $N(C_1)$ is a minimal separating set of $G$, so we may put $S = N(C_1)$ and let $G_1$ and $G_2$ be as defined before~\stepref{K9=Extremal:Gi8}, where $V(G_1) = V(C_1) \cup N(C_1)$.
Let $G_2' = G_2 - N[x]$.
Suppose $\{v_2, v_4, v_5, v_6\} \subseteq N(C_1)$.
Then by contracting $V(G_2') \cup \{v_3, v_8\}$ onto $z$, and by contracting the edges $v_0 x$, $v_1 v_7$, and $v_2 v_9$, we see $e(G[S]) + d_1 = 21$, and by contracting $C_1$ onto $v_0$, we see $d_2 \ge 3$.
But then by~\stepref{K9=Extremal:eG}, $6n - 20 = e(G) \le 6n - 22$, a contradiction.
Thus $\{v_3, v_7, v_8, v_9\} \subseteq N(C_1)$.
If $z$ is not adjacent to any vertex of $N(C_1)$, then by contracting $G_i - S$ onto $z$, we see $d_{3 - i} \ge 6$ for $i \in \{1, 2\}$.
But then from~\stepref{K9=Extremal:eG}, $e(G) \le 6n - 22$, again a contradiction.
So $z$ has at least one neighbor in $N(C_1) \setminus \{z\}$.
By symmetry, we may suppose $z v_0 \in E(G)$.
Now by contracting $V(G_2') \cup \{v_4, v_5\}$ onto $z$, and by contracting the edges $x v_0$, $v_2 v_9$, and $v_6 v_7$, we see that $e(G[S]) + d_1 \ge 20$.
Since $e(G) = 6n - 20$, it follows from~\stepref{K9=Extremal:eG} that $d_2 \le 2$, that is $z$ has at least four neighbors in $N(C_1) \setminus \{z\}$.
But this requires $z$ to be adjacent to both ends of some missing edge of $G[N(x)]$, contradicting~\sstepref{K9=Extremal:NbrComp:PNbr}.
\proofsquare
\smallskip

By~\stepref{K9=Extremal:deltaNx}, $v_0$ and $v_2$ have at least two common neighbors $w_1, w_2 \in V(G) \setminus N[x]$.
Similarly, $v_1$ and $v_3$ have at least two common neighbors $u_1, u_2 \in V(G) \setminus N[x]$.
By~\sstepref{K9=Extremal:NbrComp:PNbr}, the vertices $w_1, w_2, u_1, u_2$ are distinct.
By~\sstepref{K9=Extremal:NbrComp:P2Conn}, there exist two disjoint paths $P_1, P_2$ with one end in $\{w_1, w_2\}$ and the other end in $\{u_1, u_2\}$, and all internal vertices in $G - N[x]$.
By relabelling if necessary, we may assume $P_i$ has ends $w_i, u_i$ for $i \in \{1, 2\}$.
Furthermore, there exists a path $Q$ with one end in $V(P_1)$, the other end in $V(P_2)$, and all internal vertices in $G - N[x]$.
Now by contracting $P_1$ and $P_2$ each to a single vertex, contracting $Q$ to a single edge, and contracting the edges $v_1 u_1$ and $v_2 w_2$, we have added the edges $v_0 v_1, v_1 v_2, v_2 v_3$ to $G[N(x)]$, and it follows from Lemma~\ref{lem:Petersen} that $G \ge K_9^=$, a contradiction.
\proofsquare
\smallskip

\step{\label{K9=Extremal:Disconnected} $G - N[x]$ is disconnected for any vertex $x \in V(G)$ with $d(x) \in \{9, 10, 11\}$.}

Suppose there exists a vertex $x \in V(G)$ with $d(x) \in \{9, 10, 11\}$ for which $G - N[x]$ is not disconnected.
By~\stepref{K9=Extremal:NbrComp}, it follows that $G - N[x]$ must be the empty graph, that is $N[x] = V(G)$.
But then $G[N(x)]$ is a graph on at most 11 vertices with $\delta(G[N(x)]) \ge 8$, and so $e(G) \ge 4|G| > 5|G| - 14$.
Thus $G[N(x)] \ge K_8^=$ by Theorem~\ref{thm:K8=Extremal}, a contradiction.
\proofsquare
\smallskip

\step{\label{K9=Extremal:CompDeg} Let $x \in V(G)$ with $d(x) \in \{9, 10, 11\}$.
Then there is no component $C$ of $G - N[x]$ such that $d_G(y) \ge 12$ for every $y \in V(C)$.}

Suppose such a vertex $x$ and component $C$ exist.
Let $G_1 = G - V(C)$ and $G_2 = G[V(C) \cup N(C)]$.
Let $d_1$ be as defined before~\stepref{K9=Extremal:eGi}.
From~\stepref{K9=Extremal:eGi}, $e(G_2) \le 6|G_2| - 20 - d_1 = 6(|C| + |N(C)|) - 20 - d_1$.
By contracting the edge $xz$, where $z \in N(C)$ has minimum degree $d$ in $G[N(C)]$, we have $d_1 \ge |N(C)| - d - 1$, and hence $e(G_2) \le 6|C| + 5|N(C)| - 19 + d$.
Let $t = e(V(C), N(C))$.
Then $e(G_2) = e(C) + t + e(G[N(C)])$.
Since $2e(C) \ge 12|C| - t$ and $2e(G[N(C)]) \ge d|N(C)|$, we have $2e(G_2) \ge 12|C| + t + d|N(C)|$.
Hence $12|C| + 10|N(C)| - 38 + 2d \ge 2e(G_2) \ge 12|C| + t + d|N(C)|$,
which gives $-t \ge d(|N(C)| - 2) - 10|N(C)| + 38$.
Note that $G[N(x)]$ has minimum degree at least 6 by~\stepref{K9=Extremal:deltaNx}, and so $G[N(C)]$ has minimum degree at least $6 - (d_G(x) - |N(C)|)$.
Thus $d \ge |N(C)| + 6 - d_G(x) \ge |N(C)| - 5$.
Furthermore, by~\stepref{K9=Extremal:6Conn} and~\stepref{K9=Extremal:NbrComp}, $6 \le |N(C)| \le d_G(x) - 1 \le 10$.
It follows that $d(|N(C)| - 2) - 10|N(C)| > (|N(C)| - 5)(|N(C)| - 2) - 10|N(C)| = |N(C)|^2 - 17|N(C)| + 10 \ge - 62$.
Therefore $-t \ge - 24$.
Since $2e(C) \ge 12|C| - t$, we get $e(C) \ge 6|C| - 12$.
If $|C| \ge 9$, then by the minimality of $G$ we have $C \ge K_9^=$, a contradiction.
Therefore, $|C| < 9$.
From the inequality $6|C| - 12 \le e(C) \le \binom{|C|}{2}$, it follows that $|C| \le 2$.
Since every vertex of $C$ has degree at least 12 in $G$, and since $d_G(x) \le 11$, it follows that $|C| = 2$ and $d_G(x) = 11$.
But now $C$ is a component with $N(C) = N(x)$, contradicting~\stepref{K9=Extremal:NbrComp}.
\proofsquare
\smallskip

Now choose a vertex $x \in V(G)$ with $d(x) \in \{9, 10, 11\}$, such that $G - N[x]$ has a component $C$ of minimum order.
Then choose a vertex $y \in V(C)$ of least degree in $G$.
By~\stepref{K9=Extremal:deltaG9} and~\stepref{K9=Extremal:CompDeg}, we have $d_G(y) \in \{9, 10, 11\}$.
Let $K$ be the component of $G - N[y]$ containing $x$.
We claim that $N(K)$ contains all vertices of $N(y)$ that are not complete to all other vertices of $N(y)$.
Suppose not, and let $z \in N(y)$ such that $z$ has a nonneighbor in $N(y) \setminus \{z\}$ and $z \notin N(K)$.
If $z \in N(x)$, then $z \in N(K)$, and so $z \notin N(x)$, and thus $z \in V(C)$.
Therefore $d_G(z) \ge d_G(y)$ by the choice of $y$.
Thus $z$ has some neighbor $z' \in N(x) \cup V(C) \setminus N[y]$.
Now if $z' \notin V(K)$, then the component of $G - N[y]$ containing $z'$ would be a proper subgraph of $C$, contradicting our choice of $x$ and $C$.
Therefore $z' \in V(K)$, and thus $z \in N(K)$, a contradiction.
Thus $N(K)$ does contain all vertices of $N(y)$ which are not complete to all other vertices of $N(y)$, but this contradicts~\stepref{K9=Extremal:NbrComp}.
This contradiction completes the proof of Theorem~\ref{thm:K9=Extremal}.
\proofsquare

\section*{Acknowledgment}
Thank you to Zi-Xia Song for introducing me to this topic and for helpful guidance during the early stages of the project.

\end{document}